\documentclass[a4paper,12pt]{article}

\makeatletter
\@addtoreset{footnote}{page}
\makeatother

\usepackage{amsthm}
\usepackage{amsmath,amssymb,latexsym,amsfonts,mathrsfs}

\renewcommand{\d}{{\rm d}} 
\newcommand{\eps}{\ensuremath{\varepsilon}}
\newcommand{\e}{{\rm e}} 
\newcommand{\tend}[2]{\mathrel{\mathop{\longrightarrow}\limits^{#1}_{#2}}}
\renewcommand{\tilde}{\widetilde}

\newcommand{\up}{\ensuremath{\uparrow}}
\newcommand{\down}{\ensuremath{\downarrow}}

\newcommand{\rbra}[1]{\!\left( #1 \right)} 
\newcommand{\cbra}[1]{\!\left\{ #1 \right\}} 
\newcommand{\sbra}[1]{\!\left[ #1 \right]} 

\newcommand{\bC}{\ensuremath{\mathbb{C}}}

\newcommand{\bP}{\ensuremath{\mathbb{P}}}
\newcommand{\bR}{\ensuremath{\mathbb{R}}}

\newcommand{\cF}{\ensuremath{\mathcal{F}}}

\newcommand{\ve}{\ensuremath{\mbox{{\boldmath $e$}}}}
\newcommand{\vn}{\ensuremath{\mbox{{\boldmath $n$}}}}

\theoremstyle{plain}
\newtheorem{Thm}{Theorem}[section]
\newtheorem{Lem}[Thm]{Lemma}
\newtheorem{Prop}[Thm]{Proposition}

\theoremstyle{definition}

\newtheorem{Rem}[Thm]{Remark}
\newtheorem{Ex}[Thm]{Example}

\newcommand{\Proof}[2][Proof]{\begin{proof}[{#1}] #2 \end{proof}}

\numberwithin{equation}{section}

\makeatletter
\renewcommand\section{\@startsection {section}{1}{\z@}%
                                   {-3.5ex \@plus -1ex \@minus -.2ex}%
                                   {2.3ex \@plus.2ex}%
                                   {\normalfont\large\bf}}
\makeatother

\makeatletter
\renewcommand\subsection{\@startsection {subsection}{1}{\z@}%
                                   {-3.5ex \@plus -1ex \@minus -.2ex}%
                                   {2.3ex \@plus.2ex}%
                                   {\normalfont\normalsize\bf}}
\makeatother

\begin{document}

\begin{center}
{\Large \bf 
On $ h $-transforms of one-dimensional diffusions stopped upon hitting zero 
}
\end{center}
\begin{center}
Kouji Yano\footnote{
Graduate School of Science, Kyoto University, JAPAN }\footnote{
The research of this author is supported by KAKENHI (26800058) 
and partially by KAKENHI (24540390).} 
and Yuko Yano\footnote{
Department of Mathematics, Kyoto Sangyo University, JAPAN }\footnote{
The research of this author is supported by KAKENHI (23740073).}
\end{center}
\begin{center}
{\small Dedicated to the memory of Marc Yor}
\end{center}
\begin{center}
{\small \today}
\end{center}

\begin{abstract}
For a one-dimensional diffusion on an interval 
for which 0 is the regular-reflecting left boundary, 
three kinds of conditionings to avoid zero are studied. 
The limit processes are $ h $-transforms of the process stopped upon hitting zero, 
where $ h $'s are the ground state, the scale function, 
and the renormalized zero-resolvent. 
Several properties of the $ h $-transforms are investigated. 
\end{abstract}

\section{Introduction}

For the reflecting Brownian motion $ \{ (X_t),(\bP_x)_{x \in [0,\infty )} \} $ 
and its excursion measure $ \vn $ away from 0, it is well-known that 
$ \bP^0_x[X_t] = x $ for all $ x \ge 0 $ and all $ t>0 $, 
where $ \{ (X_t),(\bP^0_x)_{x \in [0,\infty )} \} $ denotes the process 
stopped upon hitting 0, 
and $ t \mapsto \vn[X_t] $ is constant in $ t>0 $. 
Here and throughout this paper we adopt the notation $ \mu[F] = \int F \d \mu $ 
for a measure $ \mu $ and a function $ F $. 
The process conditioned to avoid zero 
may be regarded as 
the $ h $-transform with respect to $ h(x) = x $ 
of the Brownian motion stopped upon hitting zero. 
The obtained process coincides with the 3-dimensional Bessel process 
and appears in various aspects of $ \vn $ 
(see, e.g., \cite{IW} and \cite{RY}). 

We study three analogues of conditioning to avoid zero 
for one-dimensional diffusion processes. 
Adopting the natural scale $ s(x)=x $, 
we let $ M = \{ (X_t)_{t \ge 0},(\bP)_{x \in I} \} $ 
be a $ D_m D_s $-diffusion on $ I $ 
where $ I'=[0,l') $ or $ [0,l'] $ 
and $ I = I' $ or $ I' \cup \{ l \} $; 
the choices of $ I' $ and $ I $ depend on $ m $ (see Section \ref{m}). 
We suppose that $ 0 $ for $ M $ is regular-reflecting. 
Let $ M^0 = \{ (X_t)_{t \ge 0},(\bP^0)_{x \in I} \} $ 
denote the process $ M $ stopped upon hitting zero. 
We focus on three functions which are involved in conditionings to avoid zero. 
The first one is the natural scale $ s(x)=x $. 
The second one is given as follows. 
When $ l' $ is natural, we set $ \gamma_* = 0 $ and $ h_* = s $. 
When $ l' $ is not natural, 
it was shown in \cite[Theorem 3.1]{MR0087012} that 
the $ q $-resolvent operator $ G^0_q $ on $ L^2(\d m) $ for $ M^0 $ is compact 
and is represented by the eigenfunction expansion 
$ G^0_q = \sum_n (q-\gamma_n)^{-1} f_n \otimes f_n $ 
with $ 0 \ge \gamma_1 > \gamma_2 > \cdots \down - \infty $; 
in this case we write $ \gamma_* = \gamma_1 $ and $ h_* = f_1 $. 
The obtained function $ h_* $ is the second one. 
The third one is 
\begin{align}
h_0(x) = \lim_{q \down 0} \{ r_q(0,0) - r_q(x,0) \} , 
\label{}
\end{align}
where $ r_q(x,y) $ denotes the resolvent density with respect to the speed measure. 
We will prove $ h_0 $ always exists 
and we call $ h_0 $ the {\em renormalized zero-resolvent}. 

We now state three theorems concerning conditionings of $ M $ to avoid zero. 
Their proofs will be given in Section \ref{C}. 
We write $ (\cF_t)_{t \ge 0} $ for the natural filtration. 
Let $ T_a $ denote the first hitting time of $ a $. 
The first conditioning is a slight generalization of 
a formula found in \cite[Section 2.2]{MR2295612}. 

\begin{Thm} \label{T*}
Let $ x \in I' \setminus \{ 0 \} $. 
Let $ T $ be a stopping time and $ F_T $ be a bounded $ \cF_T $-measurable functional. 
Then 
\begin{align}
\lim_{a \up \sup I} \frac{\bP_x[F_T ; T<T_a<T_0]}{\bP_x(T_a<T_0)} 
= \bP^0_x \sbra{ F_T \frac{X_T}{x} ; T<T_* } , 
\label{eq: T*cond}
\end{align}
where $ T_* = \sup_{a \in I} T_a $. 
(If $ l $ is an isolated point in $ I $, 
we understand that the symbol $ \lim_{a \up \sup I} $ 
means the evaluation at $ a=l $.) 
\end{Thm}

The second conditioning is essentially due to McKean \cite{MR0087012},\cite{MR0162282}. 

\begin{Thm} \label{h*}
Let $ x \in I' \setminus \{ 0 \} $. 
Let $ T $ be a stopping time and $ F_T $ be a bounded $ \cF_T $-measurable functional. 
Then 
\begin{align}
\lim_{t \to \infty } \frac{\bP_x[F_T ; T < t < T_0]}{\bP_x(t<T_0)} 
= \bP^0_x \sbra{ F_T \frac{\e^{- \gamma_* T} h_*(X_T)}{h_*(x)} ; T<\infty } . 
\label{eq: McKean0}
\end{align}
\end{Thm}

The third conditioning is 
an analogue of Doney \cite[Section 8]{Don} (see also Chaumont--Doney \cite{MR2164035}) 
for L\'evy processes. 
For $ q>0 $, 
we write $ \ve_q $ for the exponential variable independent of $ M $. 

\begin{Thm} \label{Csp}
Let $ x \in I' \setminus \{ 0 \} $. 
Let $ T $ be a stopping time and $ F_T $ be a bounded $ \cF_T $-measurable functional. 
Then 
\begin{align}
\lim_{q \down 0} \frac{\bP_x[F_T ; T<\ve_q<T_0]}{\bP_x(\ve_q<T_0)} 
= \bP^0_x \sbra{ F_T \frac{h_0(X_T)}{h_0(x)} ; T<\infty } . 
\label{eq: cond limit}
\end{align}
\end{Thm}

The aim of this paper is to investigate 
several properties of the three functions $ h_* $, $ h_0 $ and $ s $ 
and of the corresponding $ h $-transforms.

We summarize some properties of the $ h $-transforms of $ M^0 $ as follows 
(See Section \ref{m} for the definition of the boundary classification 
and see the end of Section \ref{h} for the classification of recurrence of 0; 
here we note that $ m(\infty )<\infty $ if and only if 0 is positive recurrent): 
\begin{enumerate}
\item 
If $ m(\infty )=\infty $, 
we have that $ s $, $ h_* $ and $ h_0 $ all coincide. 
If $ l' $ for $ M $ is natural with $ m(\infty )<\infty $, 
we have that $ s $ and $ h_* $ coincide. 
\item 
For the $ h $-transform of $ M^0 $ for $ h=s $, $ h_* $ or $ h_0 $, 
the boundary 0 is entrance. 
\item 
For the $ h $-transform of $ M^0 $ for $ h=s $, 
\begin{enumerate}
\item 
the process explodes to $ \infty $ in finite time 
when $ l' $ for $ M $ is entrance; 
\item 
the process has no killing inside the interior of $ I $ 
and is elastic at $ l' $ 
when $ l' $ for $ M $ is regular-reflecting; 
\item 
the process is conservative otherwise. 
\end{enumerate}
\item 
For the $ h $-transform of $ M^0 $ for $ h=h_* $, 
the process is conservative. 
\item 
For the $ h $-transform of $ M^0 $ for $ h=h_0 $ when $ m(\infty )<\infty $, 
the process has killing inside. 
\end{enumerate}

Let us give an example where the three functions are distinct from each other. 
Let $ M $ be a reflecting Brownian motion on $ [0,l'] $ 
where both boundaries 0 and $ l' $ are regular-reflecting. 
Then we have 
\begin{align}
h_*(x) = \frac{2 l'}{\pi} \sin \frac{\pi x}{2 l'} 
, \quad 
h_0(x) = x - \frac{x^2}{2 l'} 
, \quad x \in [0,l'] . 
\label{}
\end{align}
We shall come back to this example in Example \ref{ex: rb}. 

We give several remarks about earlier studies related to 
the $ h $-transforms for the three functions. 
\begin{enumerate}
\item[$ 1^{\circ}) $.] 
The $ h $-transform of $ M^0 $ for $ h=s $ 
is sometimes used to obtain a integral representation of the excursion measure: 
see Salminen \cite{MR757463}, Yano \cite{MR2266999} 
and Salminen--Vallois--Yor \cite{MR2295612}. 
\item[$ 2^{\circ}) $.] 
The penalization problems for one-dimensional diffusions 
which generalize Theorem \ref{h*} 
were studied in Profeta \cite{MR2760742},\cite{MR2968676}. 
\item[$ 3^{\circ}) $.] 
The counterpart of $ h_0 $ 
for one-dimensional symmetric L\'evy processes 
where every point is regular for itself 
has been introduced by Salminen--Yor \cite{MR2409011} 
who proved an analogue of the Tanaka formula. 
Yano--Yano--Yor \cite{MR2552915} and Yano \cite{MR2603019} \cite{MR2648275} 
investigated the $ h $-transform of $ M^0 $ and studied 
the penalisation problems and related problems. 
For an approach to asymmetric cases, see Yano \cite{MR3072331}. 
\end{enumerate}

This paper is organized as follows. 
We prepare notation and several basic properties 
for one-dimensional generalized diffusions in Section \ref{m} 
and for excursion measures in Section \ref{n}. 
In Section \ref{h}, we prove existence of $ h_0 $. 
Section \ref{C} 
is devoted to the proofs of Theorems \ref{T*}, \ref{h*} and \ref{Csp}. 
In Section \ref{sec: IE}, 
we study invariance and excessiveness of $ h_0 $ and $ s $. 
In Section \ref{HT}, 
we study several properties of the $ h $-transforms.

{\bf Acknowledgment.} 
The authors are thankful to Professor Masatoshi Fukushima 
for drawing their attention to the paper \cite{CF}. 
They also thank Professor Matsuyo Tomisaki 
and Dr. Christophe Profeta for their valuable comments.

\section{Notation and basic properties for generalized diffusions} \label{m}

Let $ \tilde{m} $ and $ \tilde{s} $ be strictly-increasing functions 
$ (0,l') \to \bR $ such that 
$ \tilde{m} $ is right-continuous and $ \tilde{s} $ is continuous. 
We fix a constant $ 0<c<l' $ 
(the choice of $ c $ does not affect the subsequent argument at all). 
We set 
\begin{align}
F_1 = \iint_{l'>y>x>c} \d \tilde{m}(x) \d \tilde{s}(y) 
, \quad 
F_2 = \iint_{l'>y>x>c} \d \tilde{s}(x) \d \tilde{m}(y) . 
\label{}
\end{align}
We adopt Feller's classification of the right boundary $ l' $ with a slight refinement as follows: 
\begin{enumerate}
\item \label{enu1}
If $ F_1<\infty $ and $ F_2<\infty $, then $ l' $ is called {\em regular}. 
In this case we have $ \tilde{s}(l'-)<\infty $. 
\item 
If $ F_1<\infty $ and $ F_2 = \infty $, then $ l' $ is called {\em exit}. 
In this case we have $ \tilde{s}(l'-)<\infty $. 
\item 
If $ F_1=\infty $ and $ F_2<\infty $, then $ l' $ is called {\em entrance}. 
In this case we have $ \tilde{m}(l'-)<\infty $. 
\item 
If $ F_1=\infty $ and $ F_2=\infty $, then $ l' $ is called {\em natural}. 
In this case we have either $ \tilde{s}(l'-) = \infty $ or $ \tilde{m}(l'-)=\infty $. 
There are three subcases as follows: 
\begin{enumerate}
\item 
If $ \tilde{s}(l'-)=\infty $ and $ \tilde{m}(l'-)=\infty $, 
then $ l' $ is called {\em type-1-natural}. 
\item 
If $ \tilde{s}(l'-)=\infty $ and $ \tilde{m}(l'-)<\infty $, 
then $ l' $ is called {\em type-2-natural}. 
\item 
If $ \tilde{s}(l'-)<\infty $ and $ \tilde{m}(l'-)=\infty $, 
then $ l' $ is called {\em type-3-natural} or {\em natural-approachable}. 
\end{enumerate}
\end{enumerate}
The classification of the left boundary $ 0 $ is defined in a similar way. 

Let $ m $ be a function $ [0,\infty ) \to [0,\infty ] $ 
which is non-decreasing, right-continuous and $ m(0)=0 $. 
We assume that there exist $ l' $ and $ l $ with $ 0<l' \le l \le \infty $ such that 
\begin{align}
\text{$ m $ is} \ 
\begin{cases}
\text{strictly-increasing on $ [0,l') $}, \\
\text{flat and finite on $ [l',l) $}, \\
\text{infinite on $ [l,\infty ) $}. 
\end{cases}
\label{}
\end{align}
We take $ \tilde{m} = m|_{(0,l')} $ 
and the natural scale $ \tilde{s}(x) = s(x) = x $ on $ (0,l') $ 
to adopt the classification of the boundaries 0 and $ l' $. 
We choose the intervals $ I' $ and $ I $ as follows: 
\\[-8mm]
\begin{enumerate}
\item 
If $ l' $ is regular, there are three subcases 
related to the boundary condition as follows: 
\\[-7mm]
\begin{enumerate}
\item 
If $ l'<l=\infty $, then $ l' $ is called {\em regular-reflecting} 
and $ I'=I=[0,l'] $. 
\item 
If $ l'<l<\infty $, then $ l' $ is called {\em regular-elastic}, 
$ I'=[0,l'] $ and $ I = [0,l'] \cup \{ l \} $. 
\item 
If $ l'=l<\infty $, then $ l' $ is called {\em regular-absorbing}, 
$ I'=[0,l) $ and $ I = [0,l] $. 
\\[-7mm]
\end{enumerate}
\item 
If $ l' $ is exit, then $ l'=l<\infty $, $ I'=[0,l) $ and $ I = [0,l] $. 
\item 
If $ l' $ is entrance, then $ l'=l=\infty $ and $ I'=I=[0,\infty ) $. 
\item 
If $ l' $ is natural, then $ l'=l \le \infty $ and $ I'=I=[0,l) $. 
\end{enumerate}

We always write $ (X_t)_{t \ge 0} $ for the coordinate process 
on the space of paths $ \omega : [0,\infty ) \to \bR \cup \{ \partial \} $ 
with $ \zeta(\omega) \in [0,\infty ) $ such that 
$ \omega:[0,\zeta(\omega)) \to \bR $ is continuous 
and $ \omega(t)=\partial $ for all $ t \ge \zeta(\omega) $. 
We always adopt the canonical representation for each process 
and the right-continuous filtration $ (\cF_t)_{t \ge 0} $ 
defined by $ \cF_t = \bigcap_{s>t} \sigma(X_u:u \le s) $. 

We study a $ D_m D_s $-generalized diffusion on $ I $ 
where 0 is the regular-reflecting boundary (see Watanabe \cite[Section 3]{MR0365731}). 
Such a process can be constructed from the Brownian motion via the time-change method. 
Let $ \{ (X_t)_{t \ge 0} , (\bP^B_x)_{x \in \bR} \} $ denote the Brownian motion on $ \bR $ 
and let $ \ell(t,x) $ denote its jointly-continuous local time. 
Set $ A(t) = \int_I \ell(t,x) \d m(x) $ 
and write $ A^{-1} $ for the right-continuous inverse of $ A $. 
Then the process $ \{ (X_{A^{-1}(t)})_{t \ge 0} , (\bP^B_x)_{x \in I} \} $ 
is a realization of the desired generalized diffusion. 

Let $ M = \{ (X_t)_{t \ge 0},(\bP_x)_{x \in I} \} $ 
denote the $ D_m D_s $-generalized diffusion. 
We denote the resolvent operator of $ M $ by 
\begin{align}
R_qf(x) = \bP_x \sbra{ \int_0^{\infty } \e^{-qt} f(X_t) \d t } 
, \quad q>0 . 
\label{}
\end{align}
For $ x \in I $, we write 
\begin{align}
T_x = \inf \{ t>0 : X_t = x \} . 
\label{}
\end{align}
Then, for $ a,x,b \in I $ with $ a<x<b $, we have 
\begin{align}
\bP_x(T_a>T_b) = \frac{x-a}{b-a} . 
\label{}
\end{align}
Note that, whenever $ l \in I $, 
we have $ \bP_x(T_l<\infty )=1 $ for all $ x \in I $ 
and $ l $ is a trap for $ M $.

For a function $ f:[0,l) \to \bR $, we define 
\begin{align}
Jf(x) = \int_0^x \d y \int_{(0,y]} f(z) \d m(z) . 
\label{}
\end{align}
We sometimes write $ s(x) = x $ to emphasize the natural scale. 
For $ q \in \bC $, we write $ \phi_q $ and $ \psi_q $ 
for the unique solutions of the integral equations 
\begin{align}
\phi_q = 1 + q J \phi_q 
\quad \text{and} \quad 
\psi_q = s + q J \psi_q 
\quad \text{on $ [0,l) $}, 
\label{}
\end{align}
respectively. They can be represented as 
\begin{align}
\phi_q = \sum_{n=0}^{\infty } q^n J^n 1 
\quad \text{and} \quad 
\psi_q = \sum_{n=0}^{\infty } q^n J^n s . 
\label{}
\end{align}
Let $ q>0 $. Note that $ \phi_q $ and $ \psi_q $ are non-negative increasing functions. 
Set 
\begin{align}
H(q) = \lim_{x \up l} \frac{\psi_q(x)}{\phi_q(x)} 
= \int_0^l \frac{1}{\phi_q(x)^2} \d x . 
\label{eq: H(q)}
\end{align}
Then there exist $ \sigma $-finite measures $ \sigma $ and $ \sigma^* $ on $ [0,\infty ) $ 
such that 
\begin{align}
H(q) = \int_{[0,\infty )} \frac{1}{q+\xi} \sigma(\d \xi) 
\quad \text{and} \quad 
\frac{1}{q H(q)} = \int_{[0,\infty )} \frac{1}{q+\xi} \sigma^*(\d \xi) . 
\label{}
\end{align}
Note that 
\begin{align}
l = \lim_{q \down 0} H(q) = \int_{[0,\infty )} \frac{\sigma(\d \xi)}{\xi} 
= \frac{1}{\sigma^*(\{ 0 \})} 
\in (0,\infty ] . 
\label{}
\end{align}
If we write $ m(\infty ) = \lim_{x \to \infty } m(x) $, we have 
\begin{align}
\pi_0 := \lim_{q \down 0} q H(q) = \sigma(\{ 0 \}) 
= \frac{1}{\int_{[0,\infty )} \frac{\sigma^*(\d \xi)}{\xi}} 
= \frac{1}{m(\infty )} \in [0,\infty ) . 
\label{}
\end{align}
Note that $ \pi_0=0 $ whenever $ l<\infty $. 
We define 
\begin{align}
\rho_q(x) = \phi_q(x) - \frac{1}{H(q)} \psi_q(x) . 
\label{}
\end{align}
Then the function $ \rho_q $ is a non-negative decreasing function on $ [0,l) $ 
which satisfies 
\begin{align}
\rho_q = 1 - \frac{s}{H(q)} + q J \rho_q . 
\label{eq: IErhoq}
\end{align}

We define 
\begin{align}
r_q(x,y) = r_q(y,x) = H(q) \phi_q(x) \rho_q(y) 
\quad 0 \le x \le y , \ x,y \in I' . 
\label{}
\end{align}
In particular, we have $ r_q(0,x) = r_q(x,0) = H(q) \rho_q(x) $ and $ r_q(0,0) = H(q) $. 
It is well-known (see, e.g., \cite{IM}) that 
\begin{align}
\bP_x[\e^{-qT_y}] = \frac{r_q(x,y)}{r_q(y,y)} 
, \quad x,y \in I' , \ q>0 . 
\label{eq: LT Ty}
\end{align}
In particular, we have 
\begin{align}
\rho_q(x) = \frac{r_q(x,0)}{r_q(0,0)} = \bP_x[\e^{-qT_0}] 
, \quad x \in I' , \ q>0 . 
\label{}
\end{align}
We write $ M' = \{ (X_t)_{t \ge 0},(\bP_x')_{x \in I'} \} $ 
for the process $ M $ killed upon hitting $ l $. 
We write $ R'_q $ for the resolvent operator of $ M' $. 
It is well-known (see, e.g., \cite{IM}) that $ r_q(x,y) $ is the resolvent density of $ M' $ 
with respect to $ \d m $, or in other words, 
\begin{align}
R'_qf(x) = \int_{I'} f(y) r_q(x,y) \d m(y) . 
\label{eq: resol}
\end{align}
We have the resolvent equation 
\begin{align}
\int_{I'} r_q(x,y) r_p(y,z) \d m(y) = \frac{r_p(x,z)-r_q(x,z)}{q-p} 
, \quad x,z \in I' , \ q,p>0 . 
\label{}
\end{align}

If $ l \in I $, we define 
\begin{align}
r_q(l,y) =& 0 
\quad \text{for $ y \in I' $}, 
\label{} \\
r_q(x,l) =& \frac{1}{q} - R'_q1(x) 
\quad \text{for $ x \in I' $}, 
\label{} \\
r_q(l,l) =& \frac{1}{q} , 
\label{}
\end{align}
and define a measure $ \tilde{m} $ on $ I $ by 
\begin{align}
\tilde{m}(\d y) = 1_{I'}(y) \d m(y) + \delta_l(\d y) . 
\label{}
\end{align}
We emphasize that $ r_q(x,y) $ is no longer symmetric 
when either $ x $ or $ y $ equals $ l $.

\begin{Prop}
The formulae \eqref{eq: LT Ty} and \eqref{eq: resol} extend to 
\begin{align}
\bP_x[\e^{-qT_y}] =& \frac{r_q(x,y)}{r_q(y,y)} 
, \quad x,y \in I, \ q>0 , 
\label{eq: LTTy} \\
R_qf(x) =& \int_I f(y) r_q(x,y) \tilde{m}(\d y) 
, \quad x \in I , \ q>0 . 
\label{eq: resolv}
\end{align}
\end{Prop}

\Proof{
Suppose $ l \in I $. 

First, we let $ x = l $. 
Then we have 
$ \bP_l[\e^{-qT_y}] = 0 = \frac{r_q(l,y)}{r_q(y,y)} $ for $ y \in I' $ 
and 
$ \bP_l[\e^{-qT_l}] = 1 = \frac{r_q(l,l)}{r_q(l,l)} $. 
We also have 
$ R_qf(l) = \bP_l \sbra{ \int_0^{\infty } \e^{-qt} f(X_t) \d t } 
= f(l)/q = f(l) r_q(l,l) \tilde{m}(\{ l \}) $. 
Hence we obtain \eqref{eq: LTTy} and \eqref{eq: resolv} in this case. 

Second, we assume $ x \in I' $. 
On one hand, we have 
\begin{align}
\int_0^{\infty } \e^{-qt} \bP_x(t \ge T_l) \d t 
= \frac{1}{q} \bP_x[\e^{-qT_l}] = r_q(l,l) \bP_x[\e^{-qT_l}] . 
\label{}
\end{align}
On the other hand, we have 
\begin{align}
\int_0^{\infty } \e^{-qt} \bP_x(t \ge T_l) \d t 
= \int_0^{\infty } \e^{-qt} \cbra{ 1 - \bP_x(X_t \in I') } \d t 
= \frac{1}{q} - R'_q1(x) 
= r_q(x,l) . 
\label{}
\end{align}
Hence we obtain \eqref{eq: LTTy} for $ y=l $. 
Using \eqref{eq: resol}, we obtain 
\begin{align}
R_qf(x) 
=& R'_qf(x) 
+ f(l) \int_0^{\infty } \e^{-qt} \bP_x(t \ge T_l) \d t 
\label{} \\
=& \int_{I'} f(y) r_q(x,y) \d m(y) + f(l) r_q(x,l) \tilde{m}(\{ l \}) , 
\label{}
\end{align}
which implies \eqref{eq: resolv}. 
}

\section{The excursion measure away from 0} \label{n}

For $ y \in I $, we write $ (L_t(y))_{t \ge 0} $ for the local time at $ y $ 
normalized as follows (see \cite{MR1439516}): 
\begin{align}
\bP_x \sbra{ \int_0^{\infty } \e^{-qt} \d L_t(y) } 
= r_q(x,y) 
, \quad x \in I , \ q>0 . 
\label{eq: lt}
\end{align}
We write $ L_t $ for $ L_t(0) $. 
Let $ \vn $ denote the excursion measure away from $ 0 $ 
corresponding to $ (L_t)_{t \ge 0} $ (see \cite{Blu}), 
where we adopt the convention that 
\begin{align}
X_t = 0 \ \text{for all} \ t \ge T_0 , \quad \text{$ \vn $-a.e.} 
\label{}
\end{align}
We define the functional $ N_q $ by 
\begin{align}
N_qf = \vn \sbra{ \int_0^{\infty } \e^{-qt} f(X_t) \d t } 
, \quad q>0 . 
\label{}
\end{align}
Then it is well-known (see \cite{MR701528}) 
that $ \vn $ can be characterized by the following identity: 
\begin{align}
N_qf 
= \frac{R_qf(0)}{r_q(0,0)} 
\quad \text{whenever $ f(0)=0 $}. 
\label{eq: Nqf}
\end{align}
In particular, taking $ f = 1_{I \setminus \{ 0 \}} $, we have 
\begin{align}
\vn \sbra{ 1 - \e^{-q T_0} } = \frac{1}{r_q(0,0)} = \frac{1}{H(q)} 
\label{}
\end{align}
and, by \eqref{eq: H(q)}, we have 
\begin{align}
\vn(T_0=\infty ) = \lim_{q \down 0} \frac{1}{H(q)} = \frac{1}{l} . 
\label{}
\end{align}
We write $ M^0 = \{ (X_t)_{t \ge 0},(\bP^0_x)_{x \in I} \} $ 
for the process $ M $ stopped upon hitting 0 
and write $ R^0_q $ for the resolvent operator of $ M^0 $. 
By the strong Markov property of $ M $, we have 
\begin{align}
R_qf(x) = R^0_qf(x) + \bP_x[\e^{-qT_0}] R_qf(0) . 
\label{eq: R0qf}
\end{align}
The resolvent density with respect to $ \tilde{m}(\d y) $ is given as 
\begin{align}
r^0_q(x,y) = 
r_q(x,y) - \frac{r_q(x,0)r_q(0,y)}{r_q(0,0)} 
\quad \text{for $ x,y \in I $}. 
\label{}
\end{align}
Note that $ r^0_q(x,y) = \psi_q(x) \rho_q(y) $ for $ x \le y $ and that 
\begin{align}
\bP^0_x [\e^{-q T_y}] = \frac{r^0_q(x,y)}{r^0_q(y,y)} = \frac{\psi_q(x)}{\psi_q(y)} 
\quad 
\text{for $ x,y \in I $, $ x \le y $}. 
\label{}
\end{align}
Note also that $ (L_t(y))_{t \ge 0} $ is the local time at $ y $ such that 
\begin{align}
\bP^0_x \sbra{ \int_0^{\infty } \e^{-qt} \d L_t(y) } = r^0_q(x,y) 
, \quad x,y \in I \setminus \{ 0 \} , \ q>0 . 
\label{}
\end{align}
The strong Markov property of $ \vn $ may be stated as 
\begin{align}
\vn \sbra{ F_T \, G \circ \theta_T } = \vn \sbra{ F_T \bP^0_{X_T}[G] } , 
\label{}
\end{align}
where $ T $ is a stopping time, 
$ F_T $ is a non-negative $ \cF_T $-measurable functional, 
$ G $ is a non-negative measurable functional 
such that $ 0<\vn[F_T] < \infty $ or $ 0 < \vn[G \circ \theta_T] < \infty $.

Let $ x,y \in I $ be such that $ 0 < x < y $. 
Because of the properties of excursion paths of a generalized diffusion, 
we see that $ X $ under $ \vn $ hits $ y $ if and only if 
$ X $ hits $ x $ and in addition $ X \circ \theta_{T_x} $ hits $ y $. 
Hence, by the strong Markov property of $ \vn $, we have, 
\begin{align}
\vn(T_y < \infty ) 
=& \vn \rbra{ \{ T_y < \infty  \} \circ \theta_{T_x} \cap \{ T_x < \infty \} } 
\label{} \\
=& \bP^0_x(T_y<\infty ) \vn(T_x<\infty ) 
\label{} \\
=& \bP_x(T_y<T_0) \vn(T_x<\infty ) 
\label{} \\
=& \frac{x}{y} \vn(T_x<\infty ) . 
\label{}
\end{align}
This shows that $ x \vn(T_x<\infty ) $ equals a constant $ C $ in $ x \in I \setminus \{ 0 \} $, 
so that we have 
\begin{align}
\vn(T_x<\infty ) = \frac{C}{x} 
, \quad x \in I \setminus \{ 0 \} . 
\label{}
\end{align}
If $ l \in I $, then we have 
\begin{align}
C = l \vn(T_l < \infty ) = l \vn(T_0=\infty ) = 1 . 
\label{}
\end{align}
The following theorem generalizes this fact and a result of \cite{CF}. 

\begin{Thm}[see also \cite{CF}] \label{1}
In any case, $ C=1 $. 
\end{Thm}

Theorem \ref{1} will be proved at the end of Section \ref{sec: IE}. 

The following lemma is the first step of the proof of Theorem \ref{1}.

\begin{Lem}[see also \cite{CF}] \label{2}
The constant $ C $ may be represented as 
\begin{align}
C = \lim_{t \down 0} \vn[X_t] . 
\label{eq: C formula}
\end{align}
\end{Lem}

\Proof{
By definition of $ C $, we have 
\begin{align}
C = \sup_{x \in I \setminus \{ 0 \}} x \vn(T_x<\infty ) . 
\label{}
\end{align}
Since $ \vn(t<T_x<\infty ) \up \vn(T_x<\infty ) $ as $ t \down 0 $, we have 
\begin{align}
C 
=& \sup_{x \in I \setminus \{ 0 \}} x \sup_{t>0} \vn(t < T_x < \infty ) 
\label{} \\
=& \sup_{t>0} \sup_{x \in I \setminus \{ 0 \}} x \vn(t < T_x < \infty ) 
\label{} \\
=& \lim_{t \down 0} \sup_{x \in I \setminus \{ 0 \}} x \vn(t < T_x < \infty ) . 
\label{}
\end{align}
Because of the properties of excursion paths of a generalized diffusion, 
we see that $ X $ under $ \vn $ hits $ T_x $ after $ t $ if and only if 
$ X $ does not hit $ x $ nor $ 0 $ until $ t $ 
and $ X \circ \theta_t $ hits $ x $. 
Hence, by the strong Markov property of $ \vn $, we have, 
\begin{align}
\sup_{x \in I \setminus \{ 0 \}} x \vn(t < T_x < \infty ) 
=& \sup_{x \in I \setminus \{ 0 \}} x \vn(\{ T_x < \infty \} \circ \theta_t \cap \{ t< T_x \wedge T_0 \}) 
\label{} \\
=& \sup_{x \in I \setminus \{ 0 \}} \vn \sbra{ x \bP_{X_t}(T_x < T_0) ; t< T_x \wedge T_0 } 
\label{} \\
=& \sup_{x \in I \setminus \{ 0 \}} \vn \sbra{ X_t ; t< T_x \wedge T_0 } . 
\label{eq: supXtTxT0}
\end{align}
We divide the remainder of the proof into three cases. 

(i) The case $ l<\infty $. Since $ T_x \le T_l $ for $ x \in I \setminus \{ 0 \} $, we have 
\begin{align}
\text{\eqref{eq: supXtTxT0}} 
=& \vn \sbra{ X_t ; t< T_l \wedge T_0 } 
\label{} \\
=& \vn \sbra{ X_t ; t< T_0 } 
- \vn \sbra{ X_t ; T_l \le t < T_0 } . 
\label{}
\end{align}
Since $ \vn(T_l<\infty ) < \infty $, we may apply the dominated convergence theorem to obtain 
\begin{align}
\vn \sbra{ X_t ; T_l \le t < T_0 } 
\le \vn \sbra{ X_t ; T_l < \infty } \tend{}{t \down 0} 0 , 
\label{}
\end{align}
which implies Equality \eqref{eq: C formula}, since $ \vn[X_t;t<T_0] = \vn[X_t] $. 

(ii) The case $ l'<l=\infty $. 
The proof of Case (i) works if we replace $ l $ by $ l' $. 

(iii) The case $ l'=l=\infty $. Since $ T_x \up \infty $ as $ x \to \infty $, we have 
\begin{align}
\text{\eqref{eq: supXtTxT0}} 
= \lim_{x \to \infty } \vn \sbra{ X_t ; t< T_x \wedge T_0 } 
= \vn \sbra{ X_t ; t< T_0 } 
\label{}
\end{align}
by the monotone convergence theorem. 
This implies 
Equality \eqref{eq: C formula}. 
}

\section{The renormalized zero resolvent} \label{h}

For $ q>0 $ and $ x \in I $, we set 
\begin{align}
h_q(x) = r_q(0,0) - r_q(x,0) . 
\label{}
\end{align}
Note that $ h_q(x) $ is always non-negative, since we have, by \eqref{eq: LTTy}, 
\begin{align}
\frac{h_q(x)}{H(q)} = \bP_x[1-\e^{-qT_0}] . 
\label{eq: LTTy+}
\end{align}
The following theorem asserts that 
the limit $ h_0 := \lim_{q \down 0} h_q $ exists, 
which will be called the {\em renormalized zero resolvent}. 

\begin{Thm} \label{h0}
For $ x \in I $, the limit $ h_0(x) := \lim_{q \down 0} h_q(x) $ exists and is represented as 
\begin{align}
h_0(x) = s(x) - g(x) = x - g(x) , 
\label{eq: h0}
\end{align}
where 
\begin{align}
g(x) = \pi_0 J1(x) = \pi_0 \int_0^x m(y) \d y . 
\label{eq: g(x)}
\end{align}
The function $ h_0(x) $ is continuous increasing in $ x \in I $, 
positive in $ x \in I \setminus \{ 0 \} $ and zero at $ x=0 $. 
In particular, if $ \pi_0=0 $, then $ h_0 $ coincides with the scale function, i.e., 
$ h_0(x)=s(x)=x $D
\end{Thm}

\Proof{
For $ x \in I' $, we have 
\begin{align}
h_q(x) 
= H(q) \{ 1-\rho_q(x) \} 
= x - qH(q) J\rho_q(x) 
\tend{}{q \down 0} 
x - \pi_0 J1(x) , 
\label{}
\end{align}
where we used the facts that $ 0 \le \rho_q(x) \le 1 $ 
and $ \rho_q(x) \to 1 - \frac{x}{l} $($ =1 $ if $ \pi_0>0 $) as $ q \down 0 $ 
and used the dominated convergence theorem. 
If $ l \in I $, we have 
\begin{align}
h_q(l) = r_q(0,0) = H(q) \tend{}{q \down 0} l , 
\label{}
\end{align}
and hence we obtain $ h_0(l)=l $, 
which shows \eqref{eq: h0} for $ x=l $, 
since $ \pi_0=0 $ in this case. 

It is obvious that $ h_0 $ is continuous. 
If $ \pi_0=0 $, then $ h_0(x)=x $ is increasing in $ x \in I $ 
and positive in $ x \in I \setminus \{ 0 \} $. 
If $ \pi_0>0 $, then we have $ \pi_0 m(y) \le 1 $ for all $ y \in I $ 
and $ \pi_0 m(y) < 1 $ for all $ y < l' $, 
so that $ h_0(x) $ is increasing in $ x \in I $ 
and positive in $ x \in I \setminus \{ 0 \} $. 
The proof is now complete. 
}

\begin{Ex} \label{ex: rb}
Let $ 0<l'<l=\infty $ and let $ m(x) = \min \{ x,l' \} $. 
In this case, $ M $ is a Brownian motion on $ [0,l'] $ 
where both boundaries $ 0 $ and $ l' $ are regular-reflecting. 
Then we have 
\begin{align}
h_*(x) = \frac{2 l'}{\pi} \sin \frac{\pi x}{2 l'} 
, \quad 
h_0(x) = x - \frac{x^2}{2 l'} 
, \quad x \in [0,l'] . 
\label{}
\end{align}
Note that we have $ \pi_0 = 1/m(\infty ) = 1/l' $ and 
\begin{align}
\phi_q(x) =& 
\begin{cases}
\cosh \sqrt{q} x & \text{for $ x \in [0,l'] $}, \\
\phi_q(l') + \phi_q'(l')(x-l') & \text{for $ x \in (l',\infty ) $}, 
\end{cases}
\label{} \\
\psi_q(x) =& 
\begin{cases}
\sinh \sqrt{q} x / \sqrt{q} & \text{for $ x \in [0,l'] $}, \\
\psi_q(l') + \psi_q'(l')(x-l') & \text{for $ x \in (l',\infty ) $}, 
\end{cases}
\label{} \\
H(q) =& \frac{1}{\sqrt{q} \tanh \sqrt{q} l'} . 
\label{}
\end{align}
\end{Ex}

We study recurrence and transience of 0. 

\begin{Thm} \label{pos rec}
For $ M $, the following assertions hold: 
\begin{enumerate}
\item 
0 is transient if and only if $ l < \infty $. 
In this case, it holds that 
\begin{align}
\bP_x(T_0=\infty ) = \frac{x}{l} 
\quad \text{for $ x \in I $} . 
\label{}
\end{align}
\item 
0 is positive recurrent if and only if $ \pi_0>0 $. 
In this case, it holds that 
\begin{align}
\bP_x[T_0] = \frac{h_0(x)}{\pi_0} 
\quad \text{for $ x \in I $} . 
\label{}
\end{align}
\item 
0 is null recurrent if and only if $ l=\infty $ and $ \pi_0=0 $. 
\end{enumerate}
\end{Thm}

Although this theorem seems well-known, 
we give the proof for completeness of the paper. 

\Proof{
(i) By the formula \eqref{eq: LT Ty}, we have, for $ x \in I' $, 
\begin{align}
\bP_x(T_0=\infty ) 
= \lim_{q \down 0} \bP_x[1-\e^{-qT_0}] 
= \lim_{q \down 0} \cbra{ \frac{\psi_q(x)}{H(q)} - \{ \phi_q(x)-1 \} } 
= \frac{x}{l} . 
\label{}
\end{align}
Hence 0 is transient if and only if $ l<\infty $. 
If $ x=l \in I $, it is obvious that $ \bP_l(T_0=\infty )=1 $. 
This proves the claim. 

(ii) Since $ (1-\e^{-x})/x \up 1 $ as $ x \down 0 $, 
we may apply the monotone convergence theorem to see that 
\begin{align}
\bP_x[T_0] 
= \lim_{q \down 0} \frac{1}{q} \bP_x[1-\e^{-qT_0}] 
= \lim_{q \down 0} \frac{h_q(x)}{q r_q(0,0)} 
= \frac{h_0(x)}{\pi_0} , 
\label{}
\end{align}
for $ x \in I $. This shows that 
$ \bP_x[T_0]< \infty $ if and only if $ \pi_0>0 $, 
which proves the claim. 

(iii) This is obvious by (i) and (ii). 
}

We illustrate the classification of recurrence of 0 of Theorem \ref{pos rec} as follows: 
\begin{center}
\begin{tabular}{l||l|l}
 & $ l=\infty $ & $ l<\infty $ \\
\hline \hline
$ \pi_0=0 $ & (1) null recurrent & (3) transient \\
\hline
$ \pi_0>0 $ & (2) positive recurrent & impossible 
\end{tabular}

\begin{tabular}{l}
(1) $ l' $ is type-1-natural. \\
(2) $ l' $ is type-2-natural, entrance or regular-reflecting. \\
(3) $ l' $ is type-3-natural, exit, regular-elastic or regular-absorbing. 
\end{tabular}
\end{center}

\section{Various conditionings to avoid zero} \label{C}

We prove the three theorems concerning conditionings to avoid zero. 
We need the following lemma in later use. 

\begin{Lem} \label{BGfinite}
For any stopping time $ T $ and for any $ x \in I $, it holds that 
\begin{align}
\bP^0_x[X_T ; T<\infty ] \le x . 
\label{eq: BGfinite}
\end{align}
\end{Lem}

\Proof{
By \cite[Proposition II.2.8]{BG}, it suffices to prove that 
$ \bP^0_x[X_t] \le x $ for all $ t>0 $. 

Note that $ x \le \liminf_{t \down 0} \bP^0_x[X_t] $ for all $ x \in I $ 
by Fatou's lemma. 
By the help of \cite[Corollary II.5.3]{BG}, 
it suffices to prove that 
\begin{align}
\bP^0_x[X_{T_K};T_K < \infty] \le x 
\quad \text{for $ x \in I \setminus K $} 
\label{eq: BG}
\end{align}
for all compact subset $ K $ of $ I $. 

Let $ K $ be a compact subset of $ I $ and let $ x \in I \setminus K $. 
Let $ a = \sup (K \cap (0,x)) \cup \{ 0 \} $ 
and $ b = \inf (K \cap (x,l)) \cup \{ l \} $. 
Since 0 and $ l $ are traps for $ \bP^0_x $, we have 
$ T_K = T_a \wedge T_b $ on $ \{ T_K < \infty \} $, $ \bP^0_x $-a.e. 
and thus we obtain 
\begin{align}
\bP^0_x[X_{T_K};T_K<\infty ] 
\le \bP^0_x[X_{T_a \wedge T_b}] 
= a \bP_x(T_a<T_b) + b \bP_x(T_a>T_b) 
= x , 
\label{}
\end{align}
which proves \eqref{eq: BG} for $ x \notin K $. 
Hence we obtain the desired result. 
}

First, we prove Theorem \ref{T*}. 

\Proof[Proof of Theorem \ref{T*}]{
(i) Suppose that $ l' $($ =l $) is entrance or natural. 
By the strong Markov property, we have 
\begin{align}
a \bP_x[F_T ; T<T_a<T_0] 
=& a \bP_x \sbra{ F_T \bP_{X_T}(T_a<T_0) ; T<T_a \wedge T_0 } 
\label{} \\
=& \bP_x \sbra{ F_T X_T ; T<T_a \wedge T_0 } 
\label{} \\
=& \bP^0_x \sbra{ F_T X_T ; T<T_a } 
\label{eq: condition2}
\end{align}
since $ X_T=0 $ on $ \{ T \ge T_0 \} $, $ \bP^0_x $-a.s. 
By the fact that 
$ 1_{\{ T<T_a \}} \to 1_{\{ T<\infty \}} $, $ \bP^0_x $-a.s. 
and by Lemma \ref{BGfinite}, 
we may thus apply the dominated convergence theorem 
to see that 
\eqref{eq: condition2} converges as $ a \up l $ 
to $ \bP^0_x \sbra{ F_T X_T ; T<\infty } $. 
Since $ a \bP_x(T_a<T_0) = x $, 
we obtain \eqref{eq: T*cond}. 

(ii) 
Suppose that $ l' $ is regular-elastic, regular-absorbing or exit. 
By the strong Markov property, we have 
\begin{align}
l \bP_x[F_T ; T<T_l<T_0] 
=& l \bP_x[F_T \bP_{X_T}(T_l<T_0) ; T<T_l \wedge T_0] 
\label{} \\
=& \bP^0_x[F_T X_T ; T<T_l] . 
\label{}
\end{align}
Since $ \bP_x(T_l<T_0) = x/l $, we obtain \eqref{eq: T*cond}. 

(iii) 
In the case where $ l' $ is regular-reflecting, 
the proof is the same as (ii) if we replace $ l $ by $ l' $, and so we omit it. 
}

Second, we prove Theorem \ref{h*}. 

\Proof[Proof of Theorem \ref{h*}]{
By McKean \cite{MR0087012} (see also \cite{MR2266999}), we have the following facts. 
For $ \gamma \in \bR $, let $ \psi_{\gamma} $ be the solution of the integral equation 
$ \psi_{\gamma} = s + \gamma J \psi_{\gamma} $. 
Then we have the eigendifferential expansion 
\begin{align}
r_q(x,y) = 
\int_{(-\infty ,0)} (q-\gamma)^{-1} \psi_{\gamma}(x) \psi_{\gamma}(y) \theta(\d \gamma) 
\label{}
\end{align}
for the spectral measure $ \theta $. We now have 
\begin{align}
\frac{\bP_x(T_0 \in \d t)}{\d t} = 
\int_{(-\infty ,0)} \e^{\gamma t} \psi_{\gamma}(x) \theta(\d \gamma) 
, \quad 
\frac{\vn(T_0 \in \d t)}{\d t} = 
\int_{(-\infty ,0)} \e^{\gamma t} \theta(\d \gamma) , 
\label{}
\end{align}
and, for $ r>0 $, 
\begin{align}
\lim_{t \to \infty } \frac{\bP_x(T_0>t)}{\vn(T_0>t)} = h_*(x) 
, \quad 
\lim_{t \to \infty } \frac{\vn(T_0>t-r)}{\vn(T_0>t)} = \e^{- \gamma_* r} . 
\label{eq: McKean2}
\end{align}
We note that 
$ \gamma_* $ equals the supremum of the support of $ \theta $ 
and that $ h_* = \psi_{\gamma_*} $. 
If $ l' $ is natural, exit, regular-absorbing or regular-elastic, 
we see that $ \gamma_*=0 $ and $ h_* = s $. 

By the strong Markov property, we have 
\begin{align}
\bP_x \sbra{ F_T ; T < t < T_0 } 
= \bP^0_x \sbra{ F_T \left. \bP_{X_T}(T_0>t-r) \right|_{r=T} ; T<t } . 
\label{}
\end{align}
Since we have 
\begin{align}
\vn(T_0>t) 
\ge \vn \rbra{ T_y<T_0 , \ T_0 \circ \theta_{T_y} > t } 
= \frac{1}{y} \bP_y(T_0>t) , 
\label{}
\end{align}
we have $ \bP_y(T_0>t-r) \le y \vn(T_0>t-r) $. 
Hence, by Lemma \ref{BGfinite} and by the dominated convergence theorem, we obtain 
\begin{align}
\lim_{t \to \infty } \frac{1}{\vn(T_0>t)} \bP_x \sbra{ F_T ; T < t < T_0 } 
= \bP^0_x \sbra{ F_T \e^{-\gamma_* T} h_*(X_T) ; T<\infty } . 
\label{eq: McKean1}
\end{align}
Dividing both sides of \eqref{eq: McKean1} 
by those of the first equality of \eqref{eq: McKean2}, 
we obtain \eqref{eq: McKean0}. 
}

Third, we prove Theorem \ref{Csp}. 

\Proof[Proof of Theorem \ref{Csp}]{
By \eqref{eq: LTTy+}, we have 
\begin{align}
H(q) \bP_x(\ve_q<T_0) = h_q(x) \tend{}{q \down 0} h_0(x) . 
\label{}
\end{align}
Note that 
\begin{align}
\bP_x[F_T ; T<\ve_q<T_0] 
=& \bP_x \sbra{ F_T \int_T^{\infty } 1_{\{ t<T_0 \}} q \e^{-qt} \d t } 
\label{} \\
=& \bP_x \sbra{ F_T \e^{-qT} \int_0^{\infty } 1_{\{ t+T<T_0 \}} q \e^{-qt} \d t } 
\label{} \\
=& \bP_x \sbra{ F_T \e^{-qT} ; \ve_q+T<T_0 } 
\label{} \\
=& \bP_x \sbra{ F_T \e^{-qT} 1_{\{ \ve_q<T_0 \}} \circ \theta_T ; T<T_0 } . 
\label{}
\end{align}
By the strong Markov property, we have 
\begin{align}
H(q) \bP_x[F_T ; T<\ve_q<T_0] 
=& H(q) \bP_x \sbra{ F_T \e^{-qT} \bP_{X_T}(\ve_q<T_0) ; T<T_0 } 
\label{} \\
=& \bP^0_x \sbra{ F_T \e^{-qT} h_q(X_T) ; T<\infty } , 
\label{eq: condition}
\end{align}
since $ h_q(X_T)=0 $ on $ \{ T \ge T_0 \} $, $ \bP^0_x $-a.s. 
Once the interchange of the limit and the integration is justified, 
we see that 
\eqref{eq: condition} converges as $ q \down 0 $ 
to $ \bP^0_x \sbra{ F_T h_0(X_T) ; T<\infty } $, 
and hence we obtain \eqref{eq: cond limit}. 

Let us prove $ h_q(x) \le x $ for $ q>0 $ and $ x \in I $. 
If $ x \in I' $, we use \eqref{eq: IErhoq} and we have 
\begin{align}
h_q(x) 
= H(q) \{ 1 - \rho_q(x) \} 
= x - q H(q) J \rho_q (x) 
\le x . 
\label{}
\end{align}
If $ l \in I $, we have $ h_q(l) = H(q) \le l $. 
We thus see that 
the integrand of \eqref{eq: condition} is dominated by $ X_T $. 
By Lemma \ref{BGfinite}, 
we thus see that we may apply the dominated convergence theorem, 
and therefore the proof is complete. 
}

\section{Invariance and excessiveness} \label{sec: IE}

Let us introduce notation of invariance and excessiveness. 
Let $ h $ be a non-negative measurable function on $ E $. 
\begin{enumerate}
\item 
We say $ h $ is {\em $ \alpha $-invariant} for $ M^0 $ (resp. for $ \vn $) 
($ \alpha \in \bR $) if 
$ \e^{-\alpha t} \bP^0_x[h(X_t)] = h(x) $ for all $ x \in E $ and all $ t>0 $ 
(resp. there exists a positive constant $ C $ such that 
$ \e^{-\alpha t} \vn[h(X_t)] = C $ for all $ t>0 $). 
\item 
We say $ h $ is {\em $ \alpha $-excessive} for $ M^0 $ (resp. for $ \vn $) 
($ \alpha \ge 0 $) if 
$ \e^{-\alpha t} \bP^0_x[h(X_t)] \le h(x) $ for all $ x \in E $ and all $ t>0 $ 
and $ \e^{-\alpha t} \bP^0_x[h(X_t)] \to h(x) $ as $ t \down 0 $ 
(resp. there exists a positive constant $ C $ such that 
$ \e^{-\alpha t} \vn[h(X_t)] \le C $ for all $ t>0 $ 
and $ \vn[h(X_t)] \to C $ as $ t \down 0 $). 
\item 
We say $ h $ is {\em invariant} (resp. {\em excessive}) 
when $ h $ is $ 0 $-invariant (resp. $ 0 $-excessive). 
\end{enumerate}

We give the following remarks. 
\begin{enumerate}
\item 
As a corollary of Theorem \ref{h*}, 
the function $ h_* $ is $ \gamma_* $-invariant for $ M^0 $. 
\item 
As a corollary of (i), 
the function $ s $ is invariant for $ M^0 $ 
when $ l' $ for $ M $ is natural, exit, regular-absorbing or regular-elastic. 
\item 
As a corollary of Lemma \ref{BGfinite}, 
the function $ s $ is excessive for $ M^0 $ 
when $ l' $ for $ M $ is entrance or regular-reflecting. 
\item 
As a corollary of Theorem \ref{Csp}, 
the function $ h_0 $ is excessive for $ M^0 $. 
\end{enumerate}

In this section, we prove several properties to complement these statements.

Following \cite[Section 2]{MR0092046}, we introduce the operators 
\begin{align}
D_m f(x) = \lim_{\eps,\eps' \down 0} \frac{f(x+\eps)-f(x-\eps')}{m(x+\eps)-m(x-\eps')} 
\label{}
\end{align}
whenever the limit exist. 
Note that $ f(x) = \psi_q(x) $ (resp. $ f(x) = \rho_q(x) $) 
is an increasing (resp. decreasing) 
solution of the differential equation 
$ D_m D_s f = q f $ satisfying 
$ f(0)=0 $ and $ D_s f(0) = 1 $ (resp. $ f(0)=1 $ and $ D_s f(0) = -1/H(q) $).

\begin{Thm}
The function $ h_* $ is $ \gamma_* $-invariant for $ \vn $ 
when $ l' $ for $ M $ is entrance or regular-reflecting. 
\end{Thm}

\Proof{
By \cite[Section 12]{MR0068082}), we see that 
if $ D_m D_s f = F $ and $ D_m D_s g = G $ then 
\begin{align}
D_m \{ g D_s f - f D_s g \} = g F - f G . 
\label{}
\end{align}
Hence we have 
\begin{align}
(q-\gamma_*) \psi_{\gamma_*} \rho_q 
= D_m \cbra{ \psi_{\gamma_*} D_s \rho_q - \rho_q D_s \psi_{\gamma_*} } . 
\label{}
\end{align}
Integrate both sides on $ I' $ with respect to $ \d m $, we obtain 
\begin{align}
(q-\gamma_*) \int_{I'} \psi_{\gamma_*}(x) \rho_q(x) \d m(x) = 1 . 
\label{}
\end{align}
where we used the facts that $ \rho_q(0) = \psi_{\gamma_*}'(0)=1 $, 
$ \psi_{\gamma_*}(0) = \psi_{\gamma_*}'(l') = 0 $, 
$ \psi_{\gamma_*}(l') < \infty $ and $ \rho_q'(l') = 0 $. 
This shows that 
\begin{align}
N_qh_* = \frac{R_qh_*(0)}{H(q)} 
= \int_{I'} \rho_q(x) \psi_{\gamma_*}(x) \d m(x) = \frac{1}{q-\gamma_*} . 
\label{}
\end{align}
Hence we obtain 
$ \e^{-\gamma_* t} \vn[h_*(X_t)] = 1 $ for a.e. $ t>0 $. 
For $ 0<s<t $, we see, by the $ \gamma_* $-invariance of $ h_* $ for $ M^0 $, that 
\begin{align}
\e^{-\gamma_* t} \vn[h_*(X_t)] 
= \e^{-\gamma_* t} \vn[ \bP^0_{X_s}[h_*(X_{t-s})] ] 
= \e^{-\gamma_* s} \vn[ h_*(X_s) ] , 
\label{}
\end{align}
which shows that $ t \mapsto \e^{-\gamma_* t} \vn[h_*(X_t)] $ is constant in $ t>0 $. 
Thus we obtain the desired result. 
}

For later use, we need the following lemma. 

\begin{Lem} \label{integ ineq}
For $ 0 < p < q $, it holds that 
\begin{align}
\int_{(0,l')} \rho_q(y) \psi_p(y) \d m(y) \le \frac{H(p)}{H(q)(q-p)} . 
\label{eq: lem integ ineq}
\end{align}
Consequently, it holds that $ R'_q \psi_p(x) < \infty $. 
\end{Lem}

\Proof{
Let $ x < l' $. 
Using the fact that $ \rho_p \ge 0 $ and the resolvent equation, we have 
\begin{align}
\int_{(0,x]} \rho_q(y) \psi_p(y) \d m(y) 
\le& \int_{(0,x]} \rho_q(y) H(p) \phi_p(y) \d m(y) 
\label{} \\
\le& \frac{1}{H(q) \rho_p(x)} \int_{I'} r_q(0,y) r_p(y,x) \d m(y) 
\label{} \\
=& \frac{1}{H(q) \rho_p(x)} \cdot \frac{r_p(0,x) - r_q(0,x)}{q-p} 
\label{} \\
\le& \frac{r_p(0,x)}{H(q) \rho_p(x) (q-p) } 
= \frac{H(p)}{H(q)(q-p)} . 
\label{}
\end{align}
Letting $ x \up l' $, we obtain \eqref{eq: lem integ ineq}. 
}

We compute the image of the resolvent operators of $ h_0 $. 

\begin{Prop} \label{prop: Rqh0}
For $ q>0 $ and $ x \in I $, it holds that 
\begin{align}
R_qh_0(x) 
=& \frac{h_0(x)}{q} + \frac{r_q(x,0)}{q} - \frac{\pi_0}{q^2} , 
\label{eq: qRqh0} \\
R^0_qh_0(x) 
=& \frac{h_0(x)}{q} - \frac{\pi_0}{q^2} \bP_x[1-\e^{-qT_0}] , 
\label{eq: qR0qh0} \\
N_qh_0 
=& \frac{1}{q} - \frac{\pi_0}{q^2H(q)}. 
\label{eq: qNqh0}
\end{align}
\end{Prop}

\Proof{
Suppose $ x \in I' $. 
Let $ 0<p<q/2 $. 
On one hand, by the resolvent equation, we have 
\begin{align}
R_qh_p(x) 
=& r_p(0,0) \int_I r_q(x,y) \tilde{m}(\d y) 
- \int_I r_q(x,y) r_p(y,0) \tilde{m}(\d y) 
\label{} \\
=& \frac{r_p(0,0)}{q} 
- \frac{r_p(x,0) - r_q(x,0)}{q-p} 
\label{} \\
=& \frac{h_p(x)}{q-p} + \frac{r_q(x,0)}{q-p} - \frac{p H(p)}{q(q-p)} 
\label{} \\
\tend{}{p \down 0}& 
\frac{h_0(x)}{q} + \frac{r_q(x,0)}{q} - \frac{\pi_0}{q^2} . 
\label{}
\end{align}
On the other hand, for $ y \in I' $, we have 
\begin{align}
h_p(y) 
= H(p) \{ 1 - \rho_p(y) \} 
= \psi_p(y) - H(p) \{ \phi_p(y) - 1 \} 
\le \psi_{q/2}(y) . 
\label{eq: psiineq}
\end{align}
By Lemma \ref{integ ineq}, we see by the dominated convergence theorem that 
$ R_qh_p(x) \to R_qh_0(x) $ as $ p \down 0 $. 
Hence we obtain \eqref{eq: qRqh0} for $ x \in I' $. 

Suppose $ l \in I $ and $ x=l $. 
Then we have 
\begin{align}
q R_qh_0(l) 
= q h_0(l) r_q(l,l) \tilde{m}(\{ l \}) 
= h_0(l) , 
\label{}
\end{align}
which shows \eqref{eq: qRqh0} for $ x=l $, 
since $ r_q(l,0)=0 $ and $ \pi_0=0 $ in this case. 
Thus we obtain \eqref{eq: qRqh0}. 
Using \eqref{eq: R0qf}, \eqref{eq: Nqf}, \eqref{eq: qRqh0} and \eqref{eq: LTTy}, 
we immediately obtain \eqref{eq: qR0qh0} and \eqref{eq: qNqh0}. 
}

We now obtain the image of the transition operators of $ h_0 $. 

\begin{Thm} \label{invexc}
For $ t>0 $ and $ x \in I $, it holds that 
\begin{align}
\bP^0_x[h_0(X_t)] =& h_0(x) - \pi_0 \int_0^t \bP_x(s<T_0) \d s , 
\label{eq: P0xh0} \\
\vn[h_0(X_t)] =& 1 - \pi_0 \int_0^t \d s \int_{[0,\infty )} \e^{-s \xi} \sigma^*(\d \xi) . 
\label{eq: vnh0}
\end{align}
Consequently, 
for $ M^0 $ and $ \vn $, 
it holds that $ h_0 $ is invariant when $ \pi_0=0 $ 
and that $ h_0 $ is excessive but non-invariant when $ \pi_0>0 $. 
\end{Thm}

\Proof{
By \eqref{eq: qR0qh0}, we have 
\begin{align}
R^0_qh_0(x) 
= \frac{h_0(x)}{q} - \frac{\pi_0}{q} \int_0^{\infty } \e^{-qt} \bP_x(t<T_0) \d t , 
\label{}
\end{align}
which proves \eqref{eq: P0xh0} for a.e. $ t>0 $. 
By Fatou's lemma, we see that $ \bP^0_x[h_0(X_t)] \le h_0(x) $ holds 
for all $ t>0 $ and all $ x \in I $. 
For $ 0<s<t $, we have 
\begin{align}
\bP^0_x[h_0(X_t)] = \bP^0_x \sbra{ \bP^0_{X_s}[h_0(X_{t-s})] } \le \bP^0_x[h_0(X_s)] . 
\label{}
\end{align}
This shows that $ t \mapsto \bP^0_x[h_0(X_t)] $ is non-increasing. 
Since the right-hand side of \eqref{eq: P0xh0} is continuous in $ t>0 $, 
we see that \eqref{eq: P0xh0} holds for all $ t>0 $. 

By \eqref{eq: qNqh0}, we have 
\begin{align}
N_qh_0 
=& \frac{1}{q} - \frac{\pi_0}{q} \int_{[0,\infty )} \frac{1}{q+\xi} \sigma^*(\d \xi) 
\label{} \\
=& \frac{1}{q} - \frac{\pi_0}{q} \int_0^{\infty } \d t \, \e^{-qt} 
\int_{[0,\infty )} \e^{-t \xi} \sigma^*(\d \xi) , 
\label{}
\end{align}
which proves \eqref{eq: vnh0} for a.e. $ t>0 $. 
For $ 0<s<t $, we have 
\begin{align}
\vn[h_0(X_t)] = \vn \sbra{ \bP^0_{X_s}[h_0(X_{t-s})] } \le \vn[h_0(X_s)] , 
\label{}
\end{align}
from which we can conclude that \eqref{eq: vnh0} holds for all $ t>0 $. 
}

We have already proved that $ s $ is invariant for $ M^0 $ and $ \vn $ when $ \pi_0=0 $. 
We now study properties of $ s $ in the case where $ \pi_0>0 $. 
In the case $ l' $($ =\infty $) is entrance, 
the measure $ \bP_{l'} $ denotes the extension of $ M $ starting from $ l' $ 
constructed by a scale transform 
(see also Fukushima \cite[Section 6]{MR3161402}).

\begin{Thm} \label{bPs}
Suppose $ \pi_0>0 $. 
Then the following assertions hold: 
\begin{enumerate}
\item 
If $ l' $ is type-2-natural, then 
the scale function $ s(x)=x $ is invariant for $ M^0 $ and $ \vn $. 
\item 
If $ l' $ is entrance or regular-reflecting, then, for any $ q>0 $ and any $ t>0 $, 
\begin{align}
R^0_qs(x) =& \frac{x}{q} - \frac{\psi_q(x)}{q} \chi_q(l') , 
\label{eq: R0qs} \\
N_qs =& \frac{1}{q} \bP_{l'}[1-\e^{-q T_0}] , 
\label{eq: Nqs} \\
\vn[X_t] =& \bP_{l'}(t < T_0) , 
\label{eq: vnXt}
\end{align}
where 
\begin{align}
\chi_q(l') = 
\begin{cases}
\bP_{l'}[\e^{-q T_0}] & \text{if $ l' $ for $ M $ is entrance}, \\
\frac{1}{q} \cbra{ \frac{l'}{\psi_q(l')} - \rho_q(l') } 
& \text{if $ l' $ for $ M $ is regular-reflecting}. 
\end{cases}
\label{eq: bPsxzeta2}
\end{align}
Consequently, 
$ s(x)=x $ is excessive but non-invariant for $ M^0 $ and $ \vn $. 
\end{enumerate}
\end{Thm}

\Proof{
By \eqref{eq: IErhoq}, we have, for $ x \in I' $, 
\begin{align}
\rho_q'(x) = - \frac{1}{H(q)} + q \int_{(0,x]} \rho_q(y) \d m(y) . 
\label{}
\end{align}
Since $ I'=I $ when $ \pi_0>0 $, we have 
\begin{align}
\int_{I'} \rho_q(y) \d m(y) 
= \frac{1}{H(q)} R_q1(0) 
= \frac{1}{qH(q)} . 
\label{eq: l'rhoqy}
\end{align}
We write $ \rho_q(l') = \lim_{x \up l'} \rho_q(x) $. 
Recalling $ g $ is defined by \eqref{eq: g(x)} 
and using \eqref{eq: l'rhoqy}, we obtain 
\begin{align}
N_qg 
=& \pi_0 \int_0^{l'} \d x \, m(x) \int_{I' \setminus (0,x]} \rho_q(y) \d m(y) 
\label{} \\
=& - \frac{\pi_0}{q} \int_0^{l'} \d x \, m(x) \rho'_q(x) 
\label{} \\
=& - \frac{\pi_0}{q} \int_{I'} \d m(y) \int_y^{l'} \rho'_q(x) \d x 
\label{} \\
=& \frac{\pi_0}{q} \int_{I'} \d m(y) \{ \rho_q(y) - \rho_q(l') \} 
\label{} \\
=& \frac{\pi_0}{q} \cbra{ \frac{1}{qH(q)} - m(\infty ) \rho_q(l') } . 
\label{}
\end{align}

(i) If $ l' $ is type-2-natural, then, by \cite[Theorem 5.13.3]{ItoEssential}, 
we have $ \rho_q(l') = 0 $. 
By \eqref{eq: qNqh0}, we obtain $ N_qs = 1/q $. 
Since $ t \mapsto \vn[X_t] $ is non-decreasing, 
we obtain $ \vn[X_t]=1 $ for all $ t>0 $. 
We thus conclude that $ s $ is invariant for $ \vn $. 
The invariance of $ s $ for $ M^0 $ has already been remarked 
in the beginning of this section. 

(ii) 
We postpone the proof of \eqref{eq: R0qs} until the end of the proof of Theorem \ref{bPszeta}. 
Let us prove \eqref{eq: Nqs} and \eqref{eq: vnXt}. 

If $ l' $ is regular-reflecting, we have 
$ \rho_q(l') = \bP_{l'}[\e^{-q T_0}] $. 
If $ l' $ is entrance, 
then we may take limit as $ x \up l' $ and obtain 
\begin{align}
\rho_q(l') := \lim_{x \up l'} \rho_q(x) 
= \lim_{x \up l'} \bP_x[\e^{-q T_0}] 
= \bP_{l'}[\e^{-q T_0}] 
\label{eq: rhoql'}
\end{align}
(see Kent \cite[Section 6]{MR576891}). 
Since $ \pi_0 m(\infty ) = 1 $, we obtain 
\begin{align}
N_qs 
= N_qh_0 + N_qg 
= \frac{1}{q} \bP_{l'}[1-\e^{-q T_0}] 
= \int_0^{\infty } \e^{-qt} \bP_{l'}(t < T_0) \d t . 
\label{}
\end{align}
This proves \eqref{eq: Nqs} and 
$ \vn[X_t] = \bP_{l'}(t < T_0) $ for a.e. $ t>0 $. 
Since $ t \mapsto \bP_{l'}(t<T_0) $ is continuous (see Kent \cite[Section 6]{MR576891}) 
and by Lemma \ref{BGfinite}, 
we can employ the same argument as the proof of Theorem \ref{invexc}, 
and therefore we obtain \eqref{eq: vnXt}.

Suppose that $ s $ were invariant for $ M^0 $. 
Then we would see that $ \vn[X_t] = \vn[\bP^0_{X_s}[X_{t-s}]] = \vn[X_s] $ 
for $ 0 < s < t $, which would lead to the invariance of $ s $ for $ \vn $. 
This would be a contradiction. 
}

\begin{Rem}
An excessive function $ h $ is called {\em minimal} if, 
whenever $ u $ and $ v $ are excessive and $ h=u+v $, 
both $ u $ and $ v $ are proportional to $ h $. 
It is known (see Salminen \cite{MR757463}) that $ s $ is minimal. 
We do not know whether $ h_0 $ is minimal or not in the positive recurrent case. 
\end{Rem}

We now prove Theorem \ref{1}. 

\Proof[Proof of Theorem \ref{1}]{
In the case where $ \pi_0=0 $, we have $ h_0(x)=x $ by Theorem \ref{h0}. 
Hence, by Theorem \ref{invexc}, we see that 
$ \vn[X_t] \to 1 $ as $ t \down 0 $, which shows $ C=1 $ in this case. 

In the case where $ \pi_0>0 $, 
we obtain $ C=1 $ by Theorem \ref{bPs} and Lemma \ref{2}. 
The proof is therefore complete. 
}

\section{The $ h $-transforms of the stopped process} \label{HT}

We study $ h $-transforms in this section. 
For a measure $ \mu $ and a function $ f $, 
we write $ f \mu $ for the measure defined by $ f \mu(A) = \int_A f \d \mu $. 

Since $ h_* $ is $ \gamma_* $-invariant, 
there exists a conservative strong Markov process 
$ M^{h_*} = \{ (X_t)_{t \ge 0},(\bP^{h_*}_x)_{x \in I} \} $ such that 
\begin{align}
\bP^{h_*}_x =& \frac{\e^{- \gamma_* t} h_*(X_t)}{h_*(x)} \bP^0_x 
\quad \text{on $ \cF_t $ for $ t>0 $ and $ x \in I \setminus \{ 0 \} $} , 
\label{eq: htran1} \\
\bP^{h_*}_0 =& \e^{- \gamma_* t} h(X_t) \vn 
\quad \text{on $ \cF_t $ for $ t>0 $} . 
\label{eq: htran2}
\end{align}
We set 
\begin{align}
m^{h_*}(x) = \int_{(0,x]} h_*(y)^2 \tilde{m}(\d y) 
, \quad 
s^{h_*}(x) = \int_c^x \frac{\d y}{h_*(y)^2} , 
\label{}
\end{align}
where $ 0<c<l' $ is a fixed constant, 
We define, for $ q>0 $, 
\begin{align}
r^{h_*}_q(x,y) = 
\begin{cases}
\displaystyle 
\frac{r^0_{q+\gamma_*}(x,y)}{h(x) h(y)} 
& \text{for $ x,y \in I \setminus \{ 0 \} $}, \\
\displaystyle 
\frac{r_{q+\gamma_*}(0,y)}{h(y) r_{q+\gamma_*}(0,0)} 
& \text{for $ x=0 $ and $ y \in I \setminus \{ 0 \} $}. 
\end{cases}
\label{}
\end{align}
Then, 
we see that $ r^{h_*}_q(x,y) $ is a density of the resolvent $ R^{h_*}_q $ for $ M^{h_*} $. 

\begin{Thm} \label{Mh*}
For $ M^{h_*} $, the following assertions hold: 
\begin{enumerate}
\item 
For $ q>0 $, 
$ \phi^{h_*}_q = \frac{\psi_{q+\gamma_*}}{h_*} $ 
(resp. $ \rho^{h_*}_q = \frac{\rho_{q+\gamma_*}}{h_*} $) 
is an increasing (resp. decreasing) solution of $ D_{m^{h_*}} D_{s^{h_*}} f = qf $ 
satisfying $ f(0)=1 $ and $ D_{s^{h_*}}f(0)=0 $ 
(resp. $ f(0)=\infty $ and $ D_{s^{h_*}}f(0)=-1 $). 
\item 
$ M^{h_*} $ is the $ D_{m^{h_*}} D_{s^{h_*}} $-diffusion. 
\item 
0 for $ M^{h_*} $ is entrance. 
\item 
$ l' $ for $ M^{h_*} $ is entrance 
when $ l' $ for $ M $ is entrance; 
\\ 
$ l' $ for $ M^{h_*} $ is regular-reflecting 
when $ l' $ for $ M $ is regular-reflecting. 
\end{enumerate}
\end{Thm}

\Proof{
(i) For $ q \in \bR $ and for any function $ h $ such that $ D_m D_s h $ exists, 
we see that 
\begin{align}
D_{m^h} D_{s^h} \rbra{ \frac{\psi_{q+\alpha }}{h} } 
=& \frac{1}{h^2} D_m \cbra{ h^2 D_s \rbra{ \frac{\psi_{q+\alpha }}{h} } } 
\label{} \\
=& \frac{1}{h^2} D_m \cbra{ h D_s \psi_{q+\alpha } - \psi_{q+\alpha } D_s h } 
\label{} \\
=& \rbra{ q+\alpha - \frac{D_m D_s h}{h} } \frac{\psi_{q+\alpha }}{h} . 
\label{}
\end{align}
Taking $ h=h_* $ and $ \alpha = \gamma_* $. 
we obtain $ D_{m^{h_*}} D_{s^{h_*}} \phi^{h_*}_q = q \phi^{h_*}_q $. 
In the same way we obtain 
$ D_{m^{h_*}} D_{s^{h_*}} \rho^{h_*}_q = q \rho^{h_*}_q $. 
The initial conditions can be obtained easily. 

Claims (ii) and (iii) are obvious from (i). 

(iv) 
Suppose that $ l' $ for $ M $ is entrance or regular-reflecting. 
Then $ h_* $ is bounded, so that 
$ l' $ for $ M^{h_*} $ is of the same classification as $ l' $ for $ M $. 
Since $ M^{h_*} $ is conservative, we obtain the desired result. 
}

We now develop a general theory 
for the $ h $-transform with respect to an excessive function. 
Let $ \alpha \ge 0 $ and let $ h $ be a function on $ I $ 
which is $ \alpha $-excessive for $ M^0 $ and $ \vn $ 
which is positive on $ I \setminus \{ 0 \} $. 
Then it is well-known (see, e.g., \cite[Theorem 11.9]{CW}) that 
there exists a (possibly non-conservative) strong Markov process 
$ M^h = \{ (X_t)_{t \ge 0},(\bP^h_x)_{x \in I} \} $ such that 
\begin{align}
1_{\{ t<\zeta \}} \bP^h_x =& \frac{\e^{-\alpha t} h(X_t)}{h(x)} \bP^0_x 
\quad \text{on $ \cF_t $ for $ t>0 $ and $ x \in I \setminus \{ 0 \} $} , 
\label{eq: h-tran1} \\
1_{\{ t<\zeta \}} \bP^h_0 =& \e^{-\alpha t} h(X_t) \vn 
\quad \text{on $ \cF_t $ for $ t>0 $} . 
\label{eq: h-tran2}
\end{align}
We note that $ M^h $ becomes a diffusion 
when killed upon hitting $ l $ if $ l \in I $. 
If $ \alpha \ge 0 $, we see by \cite[Theorem 11.9]{CW} that 
the identities \eqref{eq: h-tran1} and \eqref{eq: h-tran2} 
are still valid if we replace the constant time $ t $ by a stopping time $ T $ 
and restrict both sides on $ \{ T<\infty \} $. 
We set 
\begin{align}
m^h(x) = \int_{(0,x]} h(y)^2 \tilde{m}(\d y) 
, \quad 
s^h(x) = \int_c^x \frac{\d y}{h(y)^2} , 
\label{}
\end{align}
where $ 0<c<l' $ is a fixed constant, 
We define, for $ q>0 $, 
\begin{align}
r^h_q(x,y) = 
\begin{cases}
\displaystyle 
\frac{r^0_{q+\alpha }(x,y)}{h(x) h(y)} 
& \text{for $ x,y \in I \setminus \{ 0 \} $}, \\
\displaystyle 
\frac{r_{q+\alpha }(0,y)}{h(y) r_{q+\alpha }(0,0)} 
& \text{for $ x=0 $ and $ y \in I \setminus \{ 0 \} $}. 
\end{cases}
\label{}
\end{align}
Then, 
we see that $ r^h_q(x,y) $ is a density of the resolvent $ R^h_q $ for $ M^h $.

\begin{Lem}
Suppose that $ h(x) \le \psi_{q+\alpha }(x) $ for all $ q>0 $ and all $ x \in I $. 
Define $ L^h_t(y) = L_t(y)/h(y)^2 $ for $ y \in I \setminus \{ 0 \} $. 
Then the process $ (L^h_t(y))_{t \ge 0} $ is the local time at $ y $ 
for $ M^h $ normalized as 
\begin{align}
\bP^h_x \sbra{ \int_0^{\infty } \e^{-qt} \d L^h_t(y) } = r^h_q(x,y) 
, \quad x \in I , \ y \in I \setminus \{ 0 \} . 
\label{eq: LTforMh}
\end{align}
It also holds that 
\begin{align}
\bP^h_x \sbra{ \e^{-q T_y} } = \frac{r^h_q(x,y)}{r^h_q(y,y)} 
, \quad x \in I , \ y \in I \setminus \{ 0 \} . 
\label{eq: HTforMh}
\end{align}
\end{Lem}

\Proof{
Since $ \bP^h_x $ is locally absolutely continuous with respect to $ \bP^0_x $, 
we see that $ (L^h_t(y))_{t \ge 0} $ is the local time at $ y $ for $ M^h $. 
Let $ x,y \in I \setminus \{ 0 \} $. 
For $ u \ge 0 $, we note that 
$ \eta_u(y) = \inf \{ t>0 : L_t(y) > u \} $ is a stopping time 
and that $ X_{\eta_u(y)} = y $ if $ \eta_u(y)<\infty $. 
Let $ 0 = u_0 < u_1 < \ldots < u_n $. Then, by the strong Markov property, we have 
\begin{align}
\bP^h_x \sbra{ \int_{\eta_{u_{j-1}(y)}}^{\eta_{u_j}(y)} f(t) \d L^h_t(y) } 
= \frac{1}{h(x) h(y)} \bP^0_x \sbra{ \e^{- \alpha \eta_{u_j}(y)} 
\int_{\eta_{u_{j-1}(y)}}^{\eta_{u_j}(y)} f(t) \d L_t(y) } ; 
\label{eq: Lh}
\end{align}
in fact, we have \eqref{eq: Lh} 
with restriction on $ \{ \eta_{u_n}(y) \le T_{\eps x} \} $ 
and then we obtain \eqref{eq: Lh} by letting $ \eps \down 0 $. 
Hence, by the monotone convergence theorem, we obtain 
\begin{align}
\bP^h_x \sbra{ \int_0^{\infty } f(t) \d L^h_t(y) } 
= \frac{1}{h(x) h(y)} \bP^0_x \sbra{ \int_0^{\infty } \e^{-\alpha t} f(t) \d L_t(y) } . 
\label{}
\end{align}
Letting $ f(t) = \e^{-qt} $, we obtain \eqref{eq: LTforMh} for $ x \in I \setminus \{ 0 \} $. 

Let $ x=0 $ and $ y \in I \setminus \{ 0 \} $. 
For $ p>0 $, we write $ \ve_p $ for an independent exponential time of parameter $ p $. 
By the strong Markov property, we have 
\begin{align}
\bP^h_0 \sbra{ \int_{\ve_p}^{\infty } \e^{-qt} \d L^h_t(y) } 
= \bP^h_0 \sbra{ \e^{-q \ve_p} \, r^h_q(X_{\ve_p},y) } . 
\label{eq: Ph0eqep}
\end{align}
On one hand, we have 
\begin{align}
\text{\eqref{eq: Ph0eqep}} 
\le& \bP^h_0 \sbra{ r^h_q(X_{\ve_p},y) } 
= p \int_I r^h_q(0,x) r^h_q(x,y) m^h(\d x) 
\label{} \\
=& \frac{p}{p-q} \cdot \cbra{ r^h_q(0,y) - r^h_p(0,y) } 
\tend{}{p \to \infty } r^h_q(0,y) . 
\label{}
\end{align}
On the other hand, since we have $ h(x) \le \psi_{q+\alpha }(x) $, we have 
\begin{align}
\text{\eqref{eq: Ph0eqep}} 
\ge& \bP^h_0 \sbra{ \e^{-q \ve_p} ; \ve_p < T_y } \frac{\rho_{q+\alpha }(y)}{h(y)} 
\tend{}{p \to \infty } r^h_q(0,y) . 
\label{}
\end{align}
By the monotone convergence theorem, we obtain \eqref{eq: LTforMh} for $ x=0 $. 

Using \eqref{eq: LTforMh} and using the strong Markov property, we obtain 
\begin{align}
\bP^h_x \sbra{ \e^{-q T_y} } 
= \frac{\bP^h_x \sbra{ \int_0^{\infty } \e^{-qt} \d L^h_t(y) }}
{\bP^h_y \sbra{ \int_0^{\infty } \e^{-qt} \d L^h_t(y) } } 
= \frac{r^h_q(x,y)}{r^h_q(y,y)} . 
\end{align}
This shows \eqref{eq: HTforMh}. 
}

\begin{Thm} \label{Ms}
For $ M^s $, i.e., the $ h $-transform for $ h=s $, 
the following assertions hold: 
\begin{enumerate}
\item 
For $ q>0 $, 
$ \phi^s_q = \frac{\psi_q}{s} $ (resp. $ \rho^s_q = \frac{\rho_q}{s} $) 
is an increasing (resp. decreasing) solution of $ D_{m^s} D_{s^s} f = qf $ 
satisfying $ f(0)=1 $ and $ D_{s^s}f(0)=0 $ (resp. $ f(0)=\infty $ and $ D_{s^s}f(0)=-1 $). 
\item 
$ M^s $ is the $ D_{m^s} D_{s^s} $-diffusion. 
\item 
0 for $ M^s $ is entrance. 
\item 
$ l' $ for $ M^s $ is of the same classification as $ l' $ for $ M $ 
when $ l'<\infty $, i.e., $ l' $ for $ M $ is exit, regular-absorbing, regular-elastic 
or type-3-natural; 
\\
$ l' $ for $ M^s $ is type-3-natural when $ l' $ for $ M $ is natural; 
\\
$ l' $ for $ M^s $ is exit 
when $ l' $($ =\infty $) for $ M $ is entrance with $ \int_c^{\infty } x^2 \d m(x) = \infty $; 
\\
$ l' $ for $ M^s $ is regular-absorbing 
when $ l' $($ =\infty $) for $ M $ is entrance with $ \int_c^{\infty } x^2 \d m(x) < \infty $; 
\\
$ l' $ for $ M^s $ is regular-elastic when $ l' $ for $ M $ is regular-reflecting. 
\end{enumerate}
\end{Thm}

\Proof{
Claim (i) can be obtained in the same way as the proof of (i) of Theorem \ref{Mh*}. 
Claims (ii) and (iii) are obvious from (i). 

(iv) 
Suppose $ l' $ for $ M $ is exit, regular-absorbing or regular-elastic. 
Then we have $ l'<\infty $, and hence it is obvious that 
$ l' $ for $ M^s $ is of the same classification as $ l' $ for $ M $. 

Suppose $ l' $ for $ M $ is natural. 
Then we have 
\begin{align}
\iint_{l' >y>x>c} \d m^s(x) \d s^s(y) 
= \int_{l' >x>c} x \d m(x) 
\ge \iint_{l' >y>x>c} \d x \d m(y) = \infty 
\label{}
\end{align}
and 
\begin{align}
\iint_{l' >y>x>c} \d s^s(x) \d m^s(y) 
= \int_{l' >y>c} \rbra{ \frac{1}{c} - \frac{1}{y} } y^2 \d m (y) 
\ge \int_{l' >y>2c} y \d m (y) = \infty . 
\label{}
\end{align}
Thus we see that $ l' $ for $ M^s $ is natural. 
Since $ s^s(l') = 1/c - 1/l' < \infty $, we see that 
$ l' $ for $ M^s $ is type-3-natural. 

Suppose $ l' $($ = \infty $) for $ M $ is entrance. 
Then we have 
\begin{align}
\iint_{\infty >y>x>c} \d m^s(x) \d s^s(y) 
= \iint_{\infty >y>x>c} \d x \d m(y) + c \{ m(\infty )-m(c) \} < \infty . 
\label{}
\end{align}
In addition, we have 
\begin{align}
\iint_{\infty >y>x>c} \d s^s(x) \d m^s(y) 
= \int_{\infty >y>c} \rbra{ \frac{1}{c} - \frac{1}{y} } y^2 \d m (y) , 
\label{}
\end{align}
which is finite if and only if $ \int_c^{\infty } x^2 \d m(x) $ is finite. 

Suppose $ l' $ for $ M $ is regular-reflecting. 
Then it is obvious that $ l' $ for $ M^s $ is regular. 
Since $ M^s $ has no killing inside $ [0,l') $ and 
since $ M^s $ is not conservative, 
we see that $ M^s $ has killing at $ l' $. 
Since we have 
\begin{align}
\bP^s_{l'}(T_x<\zeta) = \frac{x}{l'} \bP^0_{l'}(T_x<T_0) = \frac{x}{l'} < 1 
\quad \text{for all $ x<l' $}, 
\label{}
\end{align}
we see that $ M^s $ has killing at $ l' $. 
Thus we see that $ l' $ for $ M^s $ is regular-elastic. 
}

\begin{Rem}
When $ l'=\infty $ and $ \int_{\infty >x>c} x^2 \d m(x) < \infty $, 
the left boundary $ \infty $ is called of {\em limit circle type}. 
See Kotani \cite{MR2347790} for the spectral analysis 
involving Herglotz functions. 
\end{Rem}

\begin{Thm} \label{bPszeta}
Suppose $ l' $ for $ M $ is entrance or regular-reflecting. 
For $ M^s $, it holds that 
\begin{align}
\bP^s_x[\e^{- q \zeta}] = \frac{\psi_q(x)}{x} \chi_q(l') 
, \quad q>0 , \ x \in I' \setminus \{ 0 \} , 
\label{eq: bPsxzeta}
\end{align}
where $ \chi_q(l') $ is given by \eqref{eq: bPsxzeta2}. 
\end{Thm}

\Proof{
Suppose $ l' $ is entrance. 
Then we have 
\begin{align}
\bP^s_x[\e^{-q \zeta}] 
= \lim_{y \up l'} \bP^s_x[\e^{-q T_y}] 
= \lim_{y \up l'} \frac{y}{x} \cdot \frac{\psi_q(x)}{\psi_q(y)} . 
\label{}
\end{align}
By \cite[Theorem 5.13.3]{ItoEssential}, we have 
\begin{align}
\lim_{y \up l'} \frac{y}{\psi_q(y)} 
= \lim_{y \up l'} \frac{1}{\psi_q'(y)} 
= \rho_q(l') 
= \bP_{l'}[\e^{-qT_0}] . 
\label{}
\end{align}

Suppose $ l' $($ =\infty $) is regular-reflecting. 
Then we have 
\begin{align}
\bP^s_x[\e^{-q \zeta}] 
=& \bP^s_x[\e^{-q T_{l'}}] \bP^s_{l'}[\e^{-q \zeta}] 
= \frac{r^s_q(x,l')}{r^s_q(l',l')} \cdot \frac{1}{l'} R^0_qs(l') 
\label{} \\
=& \frac{\psi_q(x)}{x} \cdot \frac{\rho_q(l')}{\psi_q(l')} 
\cdot \int_{(0,l']} \psi_q(x) x \d m(x) . 
\label{}
\end{align}
Since $ D_m \{ \psi'_q(x) x - \psi_q(x) \} = q \psi_q(x) x $, we obtain 
\begin{align}
\int_{(0,l']} \psi_q(x) x \d m(x) 
= \frac{1}{q} \cbra{ \psi_q'(l') l' - \psi_q(l') } 
= \frac{1}{q} \cbra{ \frac{l'}{\rho_q(l')} - \psi_q(l') } . 
\label{}
\end{align}
Thus we obtain \eqref{eq: bPsxzeta}. 
}

We now give the proof of \eqref{eq: R0qs}. 

\Proof[Proof of \eqref{eq: R0qs}]{
Note that 
\begin{align}
1 - \bP^s_x[\e^{-q \zeta}] 
= q \int_0^{\infty } \d t \, \e^{-qt} \bP^s_x(\zeta > t) 
= \frac{q}{x} \int_0^{\infty } \d t \, \e^{-qt} \bP^0_x[X_t] 
= \frac{1}{x} q R^0_qs(x) . 
\label{}
\end{align}
Combining this fact with \eqref{eq: bPsxzeta}, we obtain \eqref{eq: R0qs}. 
}

\begin{Thm} \label{Mh0}
For $ M^{h_0} $, i.e., the $ h $-transform for $ h=h_0 $, 
the following assertions hold: 
\begin{enumerate}
\item 
For $ q>0 $, 
$ \phi^{h_0}_q = \frac{\psi_q}{h_0} $ (resp. $ \rho^{h_0}_q = \frac{\rho_q}{h_0} $) 
is an increasing (resp. decreasing) solution of $ D_{m^{h_0}} D_{s^{h_0}} f = qf $ 
satisfying $ f(0)=1 $ and $ D_{s^{h_0}}f(0)=0 $ 
(resp. $ f(0)=\infty $ and $ D_{s^{h_0}} f(0) = -1 $). 
\item 
$ M^{h_0} $ is the $ D_{m^{h_0}} D_{s^{h_0}} $-diffusion 
with killing measure $ \frac{\pi_0}{h_0} \d m^{h_0} $. 
\item 
0 for $ M^{h_0} $ is entrance; 
\item 
$ l' $ for $ M^{h_0} $ is natural when $ l' $ for $ M $ is type-2-natural; 
\\
$ l' $ for $ M^{h_0} $ is entrance when $ l' $ for $ M $ is entrance; 
\\
$ l' $ for $ M^{h_0} $ is regular when $ l' $ for $ M $ is regular-reflecting. 
\end{enumerate}
(For the boundary classifications for diffusions with killing measure, 
see, e.g., \cite[Chapter 4]{IM}.) 
\end{Thm}

\Proof{
Claim (i) can be obtained in the same way as the proof of (i) of Theorem \ref{Mh*}. 

(ii) For $ f = \frac{\psi_q}{h_0} $ or $ f= \frac{\rho_q}{h_0} $, we have 
\begin{align}
\rbra{ D_{m^{h_0}} D_{s^{h_0}} - \frac{\pi_0}{h_0} } f = q f , 
\label{}
\end{align}
since $ D_m D_s h_0 = - \pi_0 $. This shows (ii). 

Claim (iii) is obvious from (i). 

(iv) 
Suppose $ l' $ for $ M $ is type-2-natural. 
Then it is obvious that $ \lim_{x \up l'} \rho^{h_0}_q(x) = 0 $. 
Since we have $ D_m \{ h_0 \rho'_q - \rho_q h'_0 \} = (q h_0 + \pi_0) \rho_q $, we have 
\begin{align}
D_{s^{h_0}} \rho^{h_0}_q(x) 
= (h_0 \rho'_q - \rho_q h'_0)(x) 
= -1 + \int_{(0,x]} (q h_0(x) + \pi_0) \rho_q(x) \d m(x) . 
\label{}
\end{align}
Hence, by Proposition \ref{prop: Rqh0}, we obtain 
\begin{align}
\lim_{x \up l'} D_{s^{h_0}} \rho^{h_0}_q(x) 
= -1 + \frac{1}{H(q)} R_q(q h_0 + \pi_0)(0) 
= 0 . 
\label{}
\end{align}
Thus we see that $ l' $ for $ M^{h_0} $ is natural. 

Suppose $ l' $ for $ M $ is entrance. 
Note that 
\begin{align}
\frac{h_0(x)}{\pi_0} = x \int_{(x,\infty )} \d m(z) + \int_{(0,x]} z \d m(z) . 
\label{}
\end{align}
Since we have $ \int_{(0,\infty )} z \d m(z) < \infty $, we see that 
\begin{align}
h_0(l') := \lim_{x \up l'} h_0(x) = \pi_0 \int_{(0,\infty )} z \d m(z) < \infty . 
\label{}
\end{align}
This shows that $ l' $ for $ M^{h_0} $ is of the same classification as $ l' $ for $ M $. 

The last statement is obvious. 
}

\begin{Rem}
General discussions related to Theorems \ref{Ms} and \ref{Mh0} 
can be found in Maeno \cite{Maeno01}, Tomisaki \cite{MR2353708} and Takemura \cite{Takemura01}. 
\end{Rem}

\def\cprime{$'$} \def\cprime{$'$}

\end{document}